\documentclass{amsart}
\usepackage{amssymb}
\usepackage{amsmath}
\usepackage{graphicx}
\usepackage{color}
\usepackage{enumerate}
\newcommand{\Z}{\mathbb{Z}}
\newcommand{\R}{\mathbb{R}}
\newcommand{\crit}{\mathcal{C}}
\newcommand{\T}{\mathbb{T}}
\newcommand{\Sp}{\mathbb{S}}

\newcommand{\eps}{\varepsilon}

\newtheorem{thm}{Theorem}[section]
\newtheorem{cor}[thm]{Corollary}
\newtheorem{que}[thm]{Question}
\newtheorem{remark}{Remark}

\DeclareMathOperator{\Int}{int}
\DeclareMathOperator{\Per}{Per}
\DeclareMathOperator{\dist}{dist}

\begin{document}

\title{On dynamics of the Sierpi\'nski carpet}
\author{Jan P. Boro\'nski}
\address{National Supercomputing Center IT4Innovations, Division of the University of Ostrava,
	Institute for Research and Applications of Fuzzy Modeling,
	30. dubna 22, 70103 Ostrava,
	Czech Republic -- and -- Faculty of Applied Mathematics,
	AGH University of Science and Technology,
	al. Mickiewicza 30,
	30-059 Krak\'ow,
	Poland}
\author{Piotr Oprocha}
\email{jan.boronski@osu.cz}
\email{oprocha@agh.edu.pl}

\begin{abstract} 
We prove that the Sierpi\'nski curve admits a homeomorphism with strong mixing properties. We also prove that the constructed example does not
have Bowen's specification property.
\end{abstract}
\maketitle

\section{Introduction}
The aim of this note is to exhibit a homeomorphism of the Sierpi\'nski curve (known as the planar universal curve or Sierpi\'nski carpet) with some strong mixing properties. In 1993, Aarts and Oversteegen proved that the Sierpi\'nski curve admits a transitive homeomorphism \cite{AO}, answering a question of Gottschalk. They also showed that it does not admit a minimal one. Earlier, in 1991 Kato proved that the Sierpi\'nski curve,  does not admit expansive homeomorphisms \cite{K}. In \cite{Bis} Bi\'s, Nakayama and Walczak proved that the Sierpi\'nski curve admits a homeomorphism with positive entropy. They also showed that it admits a minimal group action (by \cite{AO} it cannot be a $\mathbb{Z}$-action). There has been quite a lot of interest in dynamical properties of the planar universal curve, also due to its occurrence as Julia sets of various complex maps (see e.g. \cite{D}). Nonetheless, we were unable to find any examples in the literature that would explicitly show homeomorphisms of the Sierp\'nski curve with chaosity beyond Devaney chaos. The writing of the note was also motivated by some recent questions. During the Workshop on Dynamical Systems and Continuum Theory at University of Vienna, in June of 2015 the following question was raised.
                                                                                                                                   
\begin{que}\label{que:1}
Suppose a 1-dimensional continuum $X$ admits a mixing homeomorphism. Must $X$ be $\frac{1}{n}$-indecomposable for some $n$?
\end{que}
Recall that a continuum $X$ is \textit{$\frac{1}{n}$-indecomposable}, if given $n$ mutually disjoint subcontinua of $X$ at least one of them must have empty interior in $X$. Note that the Sierpi\'nski curve is not $\frac{1}{n}$-indecomposable for any $n\in\mathbb{N}$. This is because it is locally connected, so every point has an arbitrarily small connected neighborhood. 

\begin{thm}\label{thm11}
The Sierpi\'nski curve $S$ admits a homeomorphism $H\colon S \to S$ such that: 
\begin{enumerate}
	\item\label{thm11:ber} $H$ has a fully supported measure $\mu$, such that $(H,\mu)$ is Bernoulli,
	\item\label{thm11:per} $H$ has a dense set of periodic points,
	\item\label{thm11:spec} $H$ does not have the specification property.
\end{enumerate}
\end{thm} 

Since every Bernoulli measure is strongly mixing, and $\mu$ in Theorem~\ref{thm11} is fully supported, we immediately obtain the following result,
answering Question~\ref{que:1} in the negative.
\begin{cor}
	The Sierpi\'nski curve $S$ admits a topologically mixing homeomorphism with dense set of periodic points.
\end{cor}
Our example is quite simple, however it relies on many nontrivial facts from topology and ergodic theory. In principle, the general strategy is very similar to the one in \cite{AO}, but the starting point is a bit different. We start with an Anosov torus diffeomorphism, which allows us to say much more about the dynamics of the constructed map. By the arguments in the proof of Theorem~\ref{thm11}, it seems very likely that Aarts-Oversteegen technique \cite{AO} which we employ here, will never lead to a map with the specification property.
This motivates the following natural question.
\begin{que}
	Does the Sierpi\'nski curve admit a homeomorphism with the specification property?
\end{que}

\section{Preliminaries}
By a dynamical system $(X,T)$ we mean a compact metric space $(X,d)$ with a continuous map $T\colon X\to X$. We identify $\T^2$ with the quotient space $\R^2/\Z^2$.

A \textit{Sierpi\'nski curve} is any set homeomorphic to 
$\Sp^2 \setminus \bigcup_{i=1}^\infty \Int D_i$ where
\begin{enumerate}[(S1)]
	\item\label{SC:1} each $D_i$ is a disc and $D_i\cap D_j=\emptyset$ for $i\neq j$,
	\item\label{SC:2} $\{D_i\}_{i=1}^\infty$ is a null sequence, i.e. the diameters of $D_i$ tend to zero, as $i\to\infty$.  
	\item\label{SC:3} $\bigcup_{i=1}^\infty D_i$ is dense in $\Sp^2$.
\end{enumerate}
Whyburn \cite{W} proved that Sierpi\'nski curve does not depend on the choice of the sequence of discs $\{D_i\}_{i=1}^\infty$,
that is any two Sierpi\'nski curves are homeomorphic.

\subsection{Topological notions of mixing}
A dynamical system $(X,T)$ is topologically mixing if for any two nonempty open sets $U,V$ there is an $N>0$
such that $T^n(U)\cap V\neq \emptyset$ for all $n\geq N$. There are many different extensions of the above property to characterize stronger mixing in the system. From the point of view of our work the following two are very important. It is not hard to see that they imply topological mixing.

In his seminal paper \cite{Bowen} Bowen introduced an important, strong version of mixing, called periodic specification property. Let $T\colon X\to X$ be a continuous onto map.
Following Bowen (cf. \cite{Denker}), we say that $(X,T)$ has the \emph{specification property} if
for any $\eps > 0$, there is a positive integer $N=N(\eps)$
such that for any integer $s\geq 2$, any 
$s$ points $y_1,\dots,y_s\in X$, and any sequence $0=j_1\leq k_1 < j_2 \leq k_2
< \dots < j_s \leq k_s$ of $2s$ integers with $j_{m+1} - k_m\geq N$ for $m= 1,\dots,s-1$, there is a point $x\in X$ such that, for each positive integer $m\leq s$ we have
$d(T^i(x),T^i(y_m))<\eps$ for all $j_m \leq i \leq k_m$.
If, in addition, we can select $x$ in such a way that $T^{k_m-j_1+N}(x)=x$
then $(X,T)$ has the \emph{periodic specification property}. Note that the problem of characterizing the relations between various types of mixing for maps in specified classes of one-dimensional continua is of high interest (e.g. see \cite{MH} and references therein).

\subsection{Invariant measures}
Let $X$ be a compact metric space with metric $d$ and let $M(X)$ be the space of
Borel probability measures on $X$ equipped with the {\it L\'evy-Prokhorov metric} $\rho$ defined by
$$\rho(\mu, \nu)=\inf\{\eps\colon \mu(A)\leq \nu(A^{\eps})+\eps \text{ for all Borel subsets } A\subset X\},$$ where $A^\delta=\{x : \dist(x,A)<\delta\}$.
While the formula defining $\rho$ is not symmetric, it is an old result of Strassen that $\rho$ is in fact a symmetric function (see Section~2.3 in \cite{Huber} and comments therein). The topology induced by $\rho$ coincides with the weak$^*$-topology on $M(X)$. It is also well known that $(M(X), \rho)$ is a compact metric space. For a dynamical system $(X,T)$ we denote by $M_T(X)$ the set of all $T$-invariant measures from $M(X)$. For more details on L\'evy-Prokhorov metric and weak*-topology the reader is referred to \cite{Huber}, and basic properties related to invariant measures (ergodicity, strong mixing, Bernoulli shift) can be found in \cite{Walters}.

\subsection{Quasi-Hyperbolic Toral Automorphisms} Let $A$ be a $2\times 2$ matrix with integer entries such that $|\det A|=1$.
Then $A^{-1}$ also has integer entries, and so $A$ induces a homeomorphism of the 2-dimensional torus $F\colon \T^2 \to \T^2$ by $F(x)=Ax\textrm{(mod 1)}$, e.g. see \cite{Brin} for more details.
Since $|\det A|=1$, every toral automorphism preserves Lebesgue measure.
It is known that the periodic points of an ergodic toral automorphism are exactly those with rational coordinates (see \cite[Proposition~24.7]{Denker}). It was first proved by Adler and Weiss for $\T^2$, and then extended by Katznelson to each $\T^n$, that if a toral automorphism is ergodic with respect to the Lebesgue measure, then it is measure-theoretically conjugate to a Bernoulli shift (e.g. see  \cite[Theorem~24.6]{Denker}).
Following Lind \cite{Lind} we say that $F$ is \emph{quasi-hyperbolic} if $A$ does not have roots of unity as eigenvalues. In dimension $2$ every quasi-hyperbolic automorphism must be hyperbolic, that is, it does not have eigenvalues on the unit circle, and has periodic specification property \cite{Lind}.

\subsection{Branched Covering from $\mathbb{T}^2$ to $\mathbb{S}^2$}\label{sec:quotient}
Take a quotient of $\mathbb{T}^2$ by the relation $J$, that identifies $(x,y)$ with $(-x,-y)$. The relation $J$ induces a branched covering map $\pi \colon \T^2\to \Sp^2$ (see e.g. \cite[p. 140]{Walters}), which is
2-to-1 except at four branch points in $\mathbb{T}^2$ given by $\crit=\{(0,0), (1/2,0), (0,1/2), (1/2,1/2)\}$.
Since the relation $J$ is preserved by any toral automorphism, for every toral automorphism $F$ we obtain a factor map $G\colon\Sp^2\to \Sp^2$
such that $G\circ \pi=\pi \circ F$. Note that if $x\not\in \crit\cup F^{-1}(\crit)$ then there is an open neighborhood $U$ of $x$ such that $U$ has at most one element of any equivalence class of the relation $J$ and the same holds for $F(U)$. Then on $U$ the factor map $\pi$ is a local isometry.

\section{Proof of Theorem~\ref{thm11}}

Start with ``Arnold's cat map'' $F\colon \mathbb{T}^2\to\mathbb{T}^2$ on the torus given by 
\begin{equation}
F(x,y)=(2x+y,x+y)\textrm{(mod 1)}.\label{ACM}
\end{equation}
Clearly $F$ is hyperbolic with eigenvalues $\lambda_1=\frac{3+\sqrt{5}}{2}$ and $\lambda_2=\frac{-\sqrt{5}-1}{2}$,
hence it has the periodic specification property.

Denote the Lebesgue measure on $\T^2$ by $\lambda$.
$F$ has a dense set of periodic points and $(\T^2,F,\lambda)$ is measure-theoretically conjugate to Bernoulli shift. Let $\pi \colon \T^2\to\Sp^2$ be the quotient map from Section~\ref{sec:quotient} and let $G$ be the induced homeomorphism of $\Sp^2$. Since $G$ is a factor of $F$, $(\Sp^2,G,\mu)$ is Bernoulli with respect to a fully supported measure $\mu$ which is a push-forward of $\lambda$ by $\pi$ (see \cite[Theorem~4.29(ii)]{Walters}). 

Fix any $x\in \T^2\setminus (\crit\cup F^{-1}(\crit))$ and identify it with $x\in [-1/2,1/2]^2$ in the universal cover. For each $v$ from the unit circle denote the corresponding radial line emerging from $x$ by $\mathcal{L}_v^x=\{x+tv : t>0\}\subset \R^2$. If we view $F$ as a linear map in the universal cover, then we clearly have that $F(\mathcal{L}_v^x)=\mathcal{L}_{w}^{F(x)}$ for $w=F(v)/||F(v)||$. Therefore for sufficiently small open neighborhoods $V,W\subset \T^2$ of $x$ and $F(x)$ there is no ambiguity in writing $F(\mathcal{L}_v^x \cap V)=\mathcal{L}_{w}^{F(x)}\cap W$. In other words $F$ locally preserves radial lines on $\T^2$ for points outside $\crit\cup F^{-1}(\crit)$. But locally near these points $\pi$ is invertible, hence also $G$ preserves radial lines locally on $\Sp^2$.

Fix a set of periodic points $\mathcal{O}\subset \Per(G)$ dense in $\Sp^2$ and such that $G(\mathcal{O})=\mathcal{O}$, $\mathcal{O}\cap \pi(\crit)=\emptyset$ and that $\Per(G)\setminus \bigcup\mathcal{O}$ is dense in $\Sp^2$. We decompose $\mathcal{O}$ into a union of $p_n$-periodic orbits $O_n$; i.e. $\mathcal{O}=\bigcup_{n\in\mathbb{N}}O_n$. 
We will modify $G$ inductively, blowing up consecutive periodic orbits from $\mathcal{O}$. Since $G$ is differentiable, this can be done by a standard procedure in differentiable dynamics (see e.g. \cite{Boyland}, p.234), or topologically adopting \cite{AO} as follows. Take the periodic obit $O_1\subset \mathcal{O}$, say $O_1=\{c,G(c),\ldots,G^{p_1-1}(c)\}$. 
Since $\pi^{-1}(O_1)\cap \crit=\emptyset$, there are open discs $D_0,\ldots, D_{p_1-1}$ such that $\pi(D_i)\cap \pi(D_j)=\emptyset$ for $i\neq j$ and $\pi$ is 1-1 on each $D_i$.
Let $U_i=\pi(D_i)$. By the discussion above we have a natural decomposition of $U_i$  into radial lines (induced locally from $D_i$) centered at $G^i(c)$, such that if $L\subset U_i$ is a sufficiently short radial line emerging from $G^i(c)$ then 
$G(L)\cap U_{i+1}$ is contained in the corresponding line. In other words, we have a decomposition of a small neighborhood of each point $G^i(c)$ into radial lines, and $G$ preserves these decompositions.

Making this formal, let $U_i$ be a small neighborhood of $G^i(c)$, and $\mathcal{F}_i=\{\mathcal{L}_v^{G^i(c)}\cap U_i: v\in \Sp^1\}$ be the family of lines emanating from the point $G^i(c)$ such that $\bigcup\mathcal{F}_i=U_i\setminus G^i(c)$. We remove $O_1$ and compactify each $\operatorname{cl} U_i\setminus \{G^{i}(c)\}$ by a topological copy $S^i_c$ of the unit circle $\mathbb{S}^1$ adding, for each index $i$ and $v\in \Sp^1\cong S^i_c$, a point $\theta_v^i \in S^i_c$ compactifying the radial line $\mathcal{L}_v^{G^i(c)}\cap U_i$. That way we obtain a $p_1$-punctured sphere $S_1$. We may easily extend $G$ to a continuous map $H_1\colon S_1\to S_1$ by setting $H_1(\theta_v^i)=\theta_{w}^j$, where $j= i+1  (\textrm{mod }p_1)$ and $w=F(v)/||F(v)||$. Clearly $H_1$ defined that way is invertible with a continuous inverse, so $H_1$ is a homeomorphism of $S_1$. Observe that the dynamics of all other points under $H_1$ in $S_1$ is exactly the same as on $\Sp^2$ for $G$. Hence we can repeat this procedure, puncturing $S_1$ and obtaining $S_2$ by replacing a periodic orbit $O_2\subset \mathcal{O}$ of length $p_2$ by a periodic sequence of circles. Proceeding inductively, we obtain a sequence of punctured spheres $S_n$ with $(\sum_{i=1}^n p_i)$-holes,
homeomorphisms $H_n\colon S_n\to S_n$ and factor maps $\pi_{n}\colon S_n \to S_{n-1}$
that collapse newly introduced circles back to points of $O_{n}$, where $S_0=\Sp^2$ and $H_0=G$.
In other words, $\pi_{n}$ reverts the modification made in step $n$. Clearly, each $\pi_n$
is a continuous onto map and $\pi_n \circ H_{n}=H_{n-1}\circ \pi_n$. By the choice of the set $\mathcal{O}$, we have $\overline{\bigcup_{n\in\mathbb{N}} O_n}=\Sp^2$.

Embed each $S_n$ in $\Sp^2$ in a natural way, and extend $\pi_n$ to a map $\eta_n\colon \Sp^2\to \Sp^2$ in the following way. If $D$ is an open disc bounded by $S_n$ in $\Sp^2$ then 
we have two possibilities. If $\pi_n(\partial D)$ is a single point then
we fix any $y\in \partial D\cap S_n$
and define $\eta_n(x)=\pi_n(y)$ for every $x\in D$. In the second case $\pi_n|_{\partial D}$ is the identity, so we can extend it to identity map on $D$. Now we shall define a 2-sphere $Q_\infty$, and a Sierpi\'nski curve $S_\infty\subseteq Q_\infty$. Denote by $S_\infty$ the inverse limit of spaces $S_n$ with bonding maps $\pi_n$ and by $Q_\infty$
the inverse limit of spheres $\Sp^2$ with $\eta_n$ as bonding maps; i.e.
\begin{eqnarray*}
Q_\infty&=&\{(z_0,z_1,\ldots):\eta_n(z_n)=z_{n-1}\},\\
S_\infty&=&\{(z_0,z_1,\ldots):\pi_n(z_n)=z_{n-1}\}\subset Q_\infty.
\end{eqnarray*}
Since each $\eta_n$ is a monotone map on a $2$-manifold, a result of Brown \cite[Theorem~4]{Brown} implies that $Q_\infty$ is homeomorphic to $\Sp^2$.
Observe that if we fix any $z\in O_n$, for some $n$, then the set $B$ of all inverse sequences in $Q_\infty$
with $z$ on the first coordinate is homeomorphic to a disc. Simply, after dropping $n$ first coordinates we see that $B$
is an inverse limit of a disk $D$ with the identity as a unique bonding map. Therefore $S_\infty$ is obtained from $Q_\infty$ by removing the interior of each element of a sequence of discs. But $\overline{\bigcup_{n\in\mathbb{N}} O_n}=\Sp^2$,
hence $S_\infty$ satisfies the conditions (S\ref{SC:1})--(S\ref{SC:3}), and so it is a Sierpi\'nski curve.
Observe that if we put $H=H_1\times H_2\times \ldots \times H_n\times \ldots$
then $H(S_\infty)=S_\infty$, therefore $H$ is a homeomorphism of the Sierpi\'nski curve.
Let $M=\Sp^2\setminus \bigcup_{n=1}^\infty O_n$ and $M_\infty=\{(z_0,z_1,\ldots)\in S_\infty : z_0\in M\}$
be the set of all inverse sequences in $S_\infty$ with the first coordinate in $M$. It follows directly from the construction
that we can view $z\in M_\infty$ as $z=(x,x,x\ldots)$ for some $x\in M$, and $H(z)=(G(x), G(x), \ldots)$.
Since periodic points of $G$ in $M$ are dense in $\Sp^2$, it is not hard to see that $H$ has a dense set of periodic points. The set $M_\infty$ is Borel, so for any Borel set $U\in S_\infty$ we can view $U\cap M_\infty$ as a Borel subset of $M$ (by projection onto the first coordinate)
and so we obtain a well defined $H$-invariant Borel probability measure $\nu$ by putting $\nu(U)=\mu(U\cap M_\infty)$. 
The measure $\mu$ is ergodic, so we have $\mu(\Sp^2\setminus M_\infty)=0$, hence also $\nu(S_\infty\setminus M_\infty)=0$, and so $\nu$ and $\mu$ are isomorphic, in particular $(S_\infty, H, \nu)$ is measure-theoretically conjugate to a Bernoulli shift. Take any open set $U$ in $S_\infty$. We claim that $\nu(U)>0$. Indeed, the basic open sets in $S_\infty$ are given by
$U_\infty=(\eta_1\circ\ldots\circ\eta_{i-1}(U_i),\ldots,\eta_{i-1}(U_i), U_i,\eta_i^{-1}(U_i),\ldots),$
where $U_i\subseteq S_i$ is open, for some $i\in\mathbb{N}$ (see e.g. Theorem 3 on p.79 in \cite{Engelking}). Since $S_i$ is a sphere with a finite number of holes, the Lebesgue measure of $U_i$ in $S_i$ is positive and $\nu(U_\infty)=\nu(U_\infty \cap M_\infty)=\mu(U_i\cap M)$, therefore $U_\infty$ has positive product measure. This shows that $\nu$ has full support, which completes the proof of Theorem~\ref{thm11}\eqref{thm11:ber}.

\begin{figure}[!ptb]
\centering
		\includegraphics[width=0.8\textwidth]{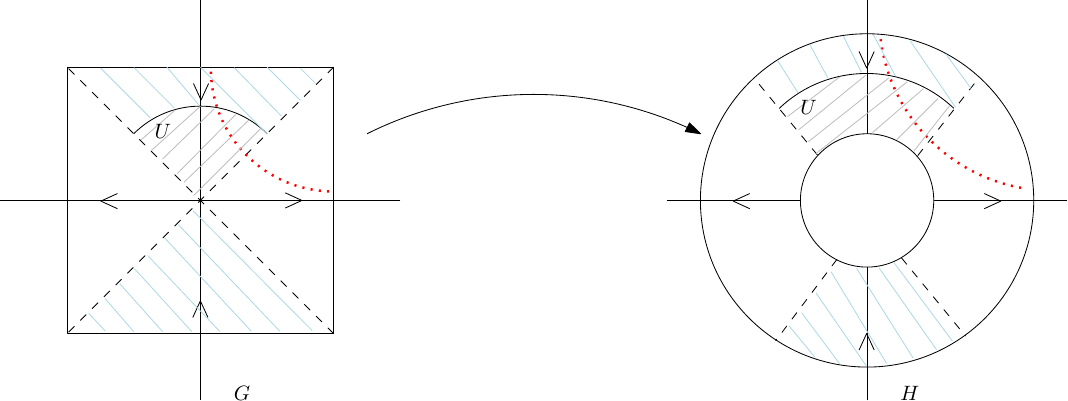}
\caption{Phase portrait for a hyperbolic point before and after a ``blow-up"}\label{fig1}
\label{phase}
\end{figure}

It remains to prove \eqref{thm11:spec}. Assume on the contrary that $H$ has the specification property. Since the specification property is preserved under higher iterations, $H^{p_1}$ has the specification property, where $p_1$ is the period of $O_1$.
For simplicity of notation replace $H$ by $H^{p_1}$ and $A$ by $A^{p_1}$.
By \cite[Theorem~2.1]{PS} the specification property implies that for every invariant measure $\mu\in M_T(S_\infty)$ there exists a sequence of ergodic measures such that $\mu_n\to \mu$,
when $n\to \infty$, in L\'evy-Prokhorov metric. Since we blew up a hyperbolic periodic point $c$ in $O_1$ in the first step of our construction, after passing from $A$ to the coordinates giving its diagonalization, we have locally a phase portrait (for $G$ and $H$) as on Figure~\ref{fig1}.
Let us start with the following observation. Consider the hyperbolic linear map $f(x,y)=(ax,by)$ where $0<a<1<b$ and $ab=1$.
Let $D=[-\eps,\eps]^2$ for some small $\eps>0$. Now let $z=(p,q)\in D$ with $|p|\geq |q|$ and assume that the trajectory of $z$ is not fully contained
in $D$. Then there exists a minimal $m\geq 1$ such that $a^m|p|\geq b^m|q|$ and $a^{m+1}|p|< b^{m+1}|q|$.
Observe that $b^{2m}|q|<\eps$ as otherwise
$\eps \leq b^{2m}|q|=a^{-m}b^m|q|\leq |p|$
and so $(p,q)\not\in D$ which is a contradiction. Now, let $v$ be a compactification of the line representing the stable direction for the hyperbolic point $c$,
and take a small neighborhood $U$ of $v$. Let $U'=\pi(U)$ where $\pi$ is the natural factor map $\pi \colon (S_\infty,H) \to (\Sp^2,G)$.
If $U$ is sufficiently small, then $\pi(U)\subset D$ and furthermore, if $(p,q)\in \pi(U)$ then $|p|\geq |q|$, see Figure~\ref{fig1}.

Fix any periodic point $u\in S_\infty$, say of period $s$, and consider the invariant measure $\hat \mu =(1-\alpha) \delta_c + \frac{\alpha}{s}\sum_{i=0}^{s-1} \delta_{H^i(u)}$
with a small $\alpha$, say $\alpha<\frac{1}{10}$. Assume also that $\pi(u)\not\in D$.
Take $0<\gamma<\alpha/2s$ small enough, so that $\dist(c, \{u,H(u),\ldots, H^{s-1}(u)\})>3\gamma $
and denote $W=B(u,2\gamma)$. We may also assume that $U$ and $\gamma$ are such small that $\overline{W}\cap \overline{U}=\emptyset$, $\pi(W)\cap D=\emptyset$, and $B(c,4\gamma)\subset U$.
Denote $V=B(c,2\gamma)$.

By results of \cite{PS} mentioned earlier, there exists an ergodic measure $\hat \nu$ such that $\rho(\hat \nu,\hat \mu)<\gamma$. 
This implies that 
$$
\hat \nu(U)\geq \hat \nu(V^\gamma)\geq \hat\mu(V)-\gamma \geq (1-\alpha)-\gamma>4/5
$$ 
and 
$$
\hat \nu(W)\geq \hat\nu (\{u\}^{2\gamma})\geq \hat \mu(\{u\}^\gamma)-\gamma>0.
$$
By the Birkhoff ergodic theorem there exists $x\in S_\infty$ such that
$\lim_{n\to \infty}\frac{1}{n}|\{j<n: H^j(x)\in U\}|=\hat \nu(U)$   and $\lim_{n\to \infty}\frac{1}{n}|\{j<n: H^j(x)\in W\}|=\hat \nu(W).$
Since $\hat \nu(W)>0$ there exists an increasing sequence $k_i$ such that $H^{k_i}(x)\in W$ for every $i$.
Let us estimate the number of iterations $k_{i}\leq j <k_{i+1}$ such that $H^j(x)\in U$.
Observe that $\pi(H^{k_i}(x))\not\in D$ and $\pi(H^{k_{i+1}}(x))\not\in D$ therefore by the earlier analysis, we see that no more than half of iterations $H^j(x)$ for $j=k_i+1,...,k_{i+1}$
can visit $U$. This implies that $\limsup_{i\to \infty} \frac{1}{k_i}|\{j<k_i:H^j(x)\in U\}|\leq 1/2$.
By the choice of $x$ we obtain that $\hat \nu(U)\leq 1/2<4/5<\hat \nu(U)$ which is a contradiction. This shows that $(S_\infty,H)$ does not have the specification property, completing the proof.
%

\begin{remark}
Similar construction works also if we start with other orientable closed surfaces, as any of them admits a branched covering onto $\mathbb{S}^2$ \cite{Alexander}. 
\end{remark}
\section{Acknowledgments}
We are grateful to an anonymous referee for careful reading, and comments that improved the paper. The authors express many thanks to L`ubomir Snoha for some helpful discussions on the properties of Besicovitch's homeomorphism used in the construction in \cite{AO}, Andrzej Bi\'s for discussion on constructions in \cite{Bis} and related topics, and Thomas Schmidt (Oregon State University) for some helpful comments. The authors' work was supported by the IT4Innovations excellence in science NPU II project LQ1602. J.~Boro\'nski also gratefully acknowledges the partial support from the MSK DT1 Support of Science and Research in the Moravian-Silesian Region  (01211/2016/RRC) ``Strengthening international cooperation in science, research and education''.
\bibliographystyle{amsplain}

\end{document}